\documentclass[12pt]{amsart}
\usepackage{amsfonts}
\usepackage{amssymb}
\theoremstyle{definition}

\theoremstyle{remark}

\newcommand{\R}{\mathbb R}

\newcommand{\ds}{\displaystyle}

\begin{document}

\centerline{\large\bf DIFFEOMORPHISMS OF PSEUDO-RIEMANNIAN MANIFOLDS}
\centerline{\large\bf AND THE VALUES OF THE CURVATURE TENSOR }
\centerline{\large\bf ON DEGENERATE PLANES
\footnote{\it SERGICA Bulgaricae mathematicae publications. Vol. 15, 1989, p. 78-86}}

\vspace{0.2in}
\centerline{\large OGNIAN T. KASSABOV}

\vspace{0.3in}
{\bf 1. Introduction.} Examining the inverse problem of the so-called "theorema egregium"
of Gauss, R. S. K\,u\,l\,k\,a\,r\,n\,i [7] has proved, that if $f$ is a 
sectional-curvature-preserving diffeomorphism of a Riemannian manifold of nowhere 
constant sectional curvature and of dimension $\ge 4$ onto another Riemannian manifold, 
then $f$ is an isometry. For other results in this direction, see [8, 9, 14].

In contrast to the definite case, when the metric is indefinite a 
sectional-curvature-preserving diffeomorphism is not necessarily an isometry. 
This occurs when the manifolds are conformal flat and of recurrent curvature on the sense 
of A. W\,a\,l\,k\,e\,r [12].

We note, that for pseudo-Riemannian manifolds the degenerate planes have an important
role. Although the sectional curvature is defined only for nondegenerate planes, 
the values of the curvature tensor on degenerate planes give a good information about 
the manifold, as it is shown e.g. in [3, 4].

In this paper it is introduced a condition, which is more general that the condition
$f$ being a sectional-curvature-preserving diffeomorphism and which is in connection
with the values of the curvature tensor on degenerate planes. It is proved, that if
this condition is fulfilled for weakly degenerate planes, then $f$ is 
sectional-curvature-preserving or the manifold is of quasi-constant curvature.
The corresponding condition for strongly degene\-rate planes implies that the diffeomorphism
is necessarily an isometry. 
 
\vspace{0.3in}
{\bf 2. Notations and preliminary results.}  Let $M$ be a (pseudo-) Riemannian
manifold with metric tensor $g$. By a plane we mean a 2-dimensional subspace of a tangent 
space. A plane $\alpha$ is said to be nondegenerate, weakly degenerate or strongly 
degenerate, if the restriction of $g$ on $\alpha$ is of rank 2, 1 or 0, respectively.
Of course, if $M$ is Riemannian, i.e. if the metric is definite, all the planes are
nondegenerate. The sectional curvature $K(\alpha)$ of a nondegenerate plane $\alpha$ 
is defined by
$$
	K(\sigma) = \frac{R(x,y,y,x)}{\pi_1(x,y,y,x)} \ ,
$$
where $\{x,y\}$ is a basis of $\alpha$, $R$ denotes the curvature tensor and $\pi_1$
is defined by
$$
	\pi_1(z,u,v,w)=g(z,w)g(u,v)-g(z,v)g(u,w) \ .
$$
 
A vetor $\xi$ in a tangent space is called isotropic, if $g(\xi,\xi)=0$ and $\xi \ne 0$. 
It is easy to see, that if a plane $\alpha$ is weakly (resp. strongly) degenerate, then there 
exists just one isotropic direction in $\alpha$ and it is orthogonal to $\alpha$ (resp.
each vector in $\alpha$ is isotropic and orthogonal to $\alpha$). On the other hand
if a plane $\alpha$ admits two orthogonal noncolinear vectors, one of which is nonisotropic
(resp. isoptropic) and the other - isotropic, then $\alpha$ is weakly (resp. strongly)
degenerate.

As usual by a curvature-like tensor on a  vector space $V$ we mean a tensor 
$T$ of type (0,4) with the properties:

1) $T(x,y,z,u)=-T(y,x,z,u)$;

2) $T(x,y,z,u)+T(y,z,x,u)+T(z,x,y,u)=0$;

3) $T(x,y,z,u)=-T(x,y,u,z)$.

The following lemma will be useful:

L\,e\,m\,m\,a A [3, 4]. {\it Let $T$ be a curvature-like tensor in a point $p$ of a 
pseudo-Riemannian manifold $M$ of dimension $\ge 3$. If\, $T$ vanishes identically on each
weakly degenerate plane in $p$, then it has the form $T=c\pi_1$ with a constant $c$.
In particular, if $T$ is the curvature tensor of $M$, then $M$ is of constant
sectional curvature in $p$.}

According to the F. Schur's theorem, if the conditions of Lemma A are fulfilled in a
connected open set $U$ with $T$ - the curvature tensor, then $U$ is of constant 
sectional curvature.

Let us recall that the Weil conformal curvature tensor $C$ for $M$ is defined by
$$
	C=R-\frac1{n-2}\varphi(S)+\frac{\tau}{(n-1)(n-2)} \pi_1 \ ,
$$
where $n={\rm dim}\, M$, $S$ is the Ricci tensor, $\tau$ is the scalar curvature and
$\varphi$ is defined by
$$
	\varphi(Q)(x,y,z,u)=g(x,u)Q(y,z)-g(x,z)Q(y,u)+g(y,z)Q(x,u)-g(y,u)Q(x,z)
$$
for any symmetric tensor $Q$ of type (0,2). It is well known [5], that if $n\ge4$
then $M$ is conformally flat if and only if $C$ vanishes identically. The following
two lemmas give criteria for vanishing of the Weil conformal curvature tensor.

L\,e\,m\,m\,a B [3, 11]. {\it Let $M$ be a (pseudo)-Riemannian manifold of dimension $\ge4$.
Assume that for every orthonormal quadruple $\{x,y,z,u\}$ of vectors in a point $p\in M$
$$
	R(x,y,z,u)=0
$$ 
holds good. Then the Weil conformal curvature tensor of $M$ vanishes in $p$. In particular,
if this holds in every point of a open set $U$, then $U$ is conformally flat.}

L\,e\,m\,m\,a C [3]. {\it Let $M$ be an $n$-dimensional pseudo-Riemannian manifold of
signature $(s,n-s)$, where $s\ge 2$, $n-s\ge 2$. If the curvature tensor of $M$ vanishes 
identically on each strongly degenerate plane in a point $p$, then the Weil conformal
curvature tensor of $M$ vanishes in $p$. In particular, if this holds in every point of 
an open set $U$, then $U$ is conformally flat.}

If $\overline M$ is another (pseudo-) Riemannian manifold, we denote the corresponding 
objects for $\overline M$ by a bar overhead. A diffeomorphism $f$ of $M$ onto $\overline M$
is said to be sectional-curvature-preserving [7], if
$$
	\overline K(f_*\alpha)=K(\alpha)
$$
for each nondegenerate plane $\alpha$ in $M$, whose image is also nondegenerate. The
corresponding condition for degenerate planes is
$$
	\lim_{\alpha\rightarrow\alpha_0} \frac{\overline K(f_*\alpha)}{K(\alpha)}=1 \ ,   \leqno (2.1)
$$
where the degenerate 2-plane \,$\alpha_0$\, is approximated by nondegenerate 2-planes, 
whose images are also nondegenerate. Note, that if $f$ is conformal, then (2.1) holds 
when and only when $f$ preserves the null sectional curvature, defined in [6].

It is easy to prove the following:

L\,e\,m\,m\,a 1. {\it Let $\alpha_0$ be a degenerate plane in $M$, such that (2.1) holds.
If $R$ doesn't vanish identically on $\alpha_0$, then $\pi_1$ vanishes identically on
$f_*\alpha_0$.}

In what follows let $M$ and $\overline M$ be pseudo-Riemannian manifolds of dimension $n$
and let $f$ be a diffeomorphism of $M$ onto $\overline M$. Then $f_*x$ will be denoted by
$\bar x$.

L\,e\,m\,m\,a 2. {\it Let $n\ge 3$ and $M$ be of nonconstant sectional curvature in a point 
$p$. If (2.1) is satisfied for each weakly degenerate plane $\alpha_0$ in the tangent
space $T_pM$, then there exists an isotropic vector $\xi$ in $p$, such that each isotropic
vector, which is sufficiently close to $\xi$ is mapped by $f_*$ in an isotropic vector in $f(p)$.}

P\,r\,o\,o\,f. According to Lemma A there exists a pair $\{x,\xi\}$ of vectors in $T_pM$
such that $x$ is an unit vector, $\xi$ is isotropic and orthogonal to $x$ and
$R(x,\xi,\xi,x)\ne 0$. Let us assume that there exists a sequence of isotropic vectors
$\xi_i$ converging to $\xi$, such that $\bar \xi_i$ are all nonisotropic. Then we
can find easily a pair  $\{x_i,\xi_i\}$ of orthogonal vectors in $T_pM$, such that
$x_i$ is a unit vector, $\xi_i$ is isotropic, $\bar\xi_i$ is nonisotropic and
$R(x_i,\xi_i,\xi_i,x_i)\ne 0$. So, without loss of generality we assume that $\bar\xi$
is nonisotropic. We may suppose also that $\bar x$ is orthogonal to $\bar\xi$. Hence,
according to Lemma 1 span\,$\{\bar x,\bar\xi\}$ is degenerate, i.e. $\bar x$ is
isotropic. Let $\xi=y+a$, where the vectors $x,\,y,\,a$ are orthogonal and 
$g(x,x)=g(y,y)=-g(a,a)=1$ or $-1$. Let $u$ be any unit vector, orthogonal to
$x,\,y,\,a$ if $n>3$ and let $u$ be the zero vector, if $n=3$. We put 
$\varepsilon=g(y,y)g(u,u)$ and $y_s=(y+su)/\sqrt{1+\varepsilon s^2}$, where $s$ is a real
number, $|s|<1$. Since $R(x,\xi,\xi,x)\ne 0$, then by Lemma 1 there exists a real
number $\delta$, $0<\delta<1$, such that for any real $s,\,t$ with $|s|,\,|t|<\delta$
we have
$$
	\bar\pi_1(\bar x+t\bar y_s, \frac{\bar y_s-t\bar x}{\sqrt{1+t^2}}+\bar a,
												\frac{\bar y_s-t\bar x}{\sqrt{1+t^2}}+\bar a,\bar x+t\bar y_s)=0 \ , 
$$    
which implies
$$
	\begin{array}{c}
			(1+t^2)\bar\pi_1(\bar x,\bar y_s,\bar y_s,\bar x)
						+2\sqrt{1+t^2}\,\bar\pi_1(\bar x,\bar y_s,\bar a,\bar x)+\bar\pi_1(\bar x,\bar a,\bar a,\bar x)  \\
		  +t^2\bar\pi_1(\bar y_s,\bar a,\bar a,\bar y_s)=0 \ .
	\end{array} \leqno (2.2)
$$
For $t=0$ this reduces to
$$
	\bar\pi_1(\bar x,\bar y_s,\bar y_s,\bar x)
						+2\bar\pi_1(\bar x,\bar y_s,\bar a,\bar x)+\bar\pi_1(\bar x,\bar a,\bar a,\bar x)  
		 =0 \ . \leqno (2.3)
$$
From (2.2) and (2.3) we find
$$
	t^2\{ \bar\pi_1(\bar x,\bar y_s,\bar y_s,\bar x) + \bar\pi_1(\bar y_s,\bar a,\bar a,\bar y_s) \}
				+2(\sqrt{1+t^2}-1)\bar\pi_1(\bar x,\bar y_s,\bar a,\bar x)=0
$$
for any real $t,\,s$ with $|t|,\,|s|<\delta$. Hence $\bar\pi_1(\bar x,\bar y_s,\bar a,\bar x)=0$.
Comparing this with (2.3) we get
$$
	\bar\pi_1(\bar x,\bar y_s,\bar y_s,\bar x)
						+\bar\pi_1(\bar x,\bar a,\bar a,\bar x)  
		 =0 
$$
and since $\bar x$ is isotropic, it is easy to conclude
$$
	\bar g(\bar x,\bar a)=0 \ ,  \leqno (2.4)
$$
$$
	\bar g(\bar x,\bar y_s)=0 \ .  \leqno (2.5)
$$
From (2.5) it follows immediately
$$
	\bar g(\bar x,\bar y)=\bar g(\bar x,\bar u)=0  \ .   \leqno (2.6)
$$
Since $\bar x$ is isotropic, (2.4) and (2.6) imply $\bar g(\bar x,\bar z)=0$ for each vector
$\bar z$ in $f(p)$ and hence $\bar x=0$, which is a contradiction. This proves the lemma.

The following assertion is an analogue of Lemma 2 for strongly degenerate planes.

L\,e\,m\,m\,a 3. {\it Let the Weil conformal curvature tensor of $M$ do not vanish
identically in $p$ and let $M$ have signature $(s,n-s)$, where $s\ge2,\,n-s\ge2$. If
$f$ satisfies (2.1) for each strongly degenerate plane $\alpha_0$ in $T_pM$, then
there exists an isotropic vector $\xi$ in $p$, such that each isotropic vector, which
is sufficiently close to $\xi$, is mapped by $f_*$ in an isotropic vector in $f(p)$.}

P\,r\,o\,o\,f. According to Lemma C there exists an orthogonal pair $\{\xi,\eta\}$
of isotropic vectors in $T_pM$, such that $R(\xi,\eta,\eta,\xi)\ne 0$. Let us assume
that there exists a sequence $\eta_i$ of isotropic vectors converging to $\eta$, such
that $\bar\eta_i$ are nonisotropic. Then we can find an orthogonal pair $\{\xi_i.\eta_i\}$
of isotropic vectors in $T_pM$, such that $\bar\eta_i$ is nonisotropic and
$R(\xi_i,\eta_i,\eta_i,\xi_i)\ne 0$. So we may asume that $\bar\eta$ is nonisotropic.
Also without loss of generality we suppose that $\xi$ and $\eta$ are orthogonal. Then Lemma 1 
implies that $\bar\xi$ is isotropic. Let $\xi=x+a,\,\eta=y+b$, where $x,\,y,\,a,\,b$ are
orthogonal and $g(x,x)=g(y,y)=-g(a,a)=-g(b,b)=1$. For any vector $z$ in $T_pM$
orthogonal to $\xi$ we put
$$
	\xi_t=\xi+t(x+z)\,;  \qquad
	\eta_t=\eta+ty \,;   \qquad
	\alpha_t={\rm span}\,\{\xi_t,\eta_t\} \ .
$$ 
Then we find
$$
	1=\lim_{t\rightarrow t_0} \frac{\overline K(\bar\alpha_t)}{K(\alpha_t)}=
	\frac{1}{R(\xi,\eta,\eta,\xi)}
	\lim_{t\rightarrow t_0}\overline R(\xi_i,\eta_i,\eta_i,\xi_i)
$$
$$
	\times\frac{t(2+t)\{2+tg(x+z,x+z)\}-t\{g(z,\eta)+tg(z,y)\}^2}
	{ \{ 2\bar g(\bar\xi,\bar x,+\bar z)+t\bar g(\bar x+\bar z,\bar x+\bar z) \} \bar g(\bar\eta_t,\bar\eta_t)
	-t\{ \bar g(\bar \xi,\bar y)+\bar g(\bar x+\bar z,\bar\eta)+t\bar g(\bar x+\bar z,\bar y)\}^2 } \ .
$$
Hence we obtain $ \bar g(\bar\xi, \bar x+\bar z)\bar g(\bar\eta,\bar\eta)=0$ or
$$
	\bar g(\bar\xi,\bar x+\bar z)=0 \ .  \leqno (2.7)
$$
In particular, if $z=0$ this implies $\bar g(\bar \xi,\bar x)=0$ and hence
$\bar g(\bar \xi,\bar a)=0$. Using again (2.7) we conclude $\bar g(\bar \xi,\bar z)=0$.
Consequently $\xi$ is orthogonal to every vector in $f(p)$, which is
a contradiction. This proves the lemma.

L\,e\,m\,m\,a 4. {\it Let $n\ge 3$ and in a point $p\in M$ there exists an isotropic 
vector $\xi$, such that each isotropic vector, which is sufficiently close to $\xi$ is 
mapped by $f_*$ into an isotropic vector. Then $f$ is a homothety in $p$.}

P\,r\,o\,o\,f. Let $\xi=x+a$, where $\{x,\,a\}$ is an orhonormal pair of vectors. Let
$y$ be an arbitrary unit vector, orthogonal to $x,\,a$. Then e.g.
$g(x,x)=g(y,y)=-g(a,a)$ and for each real $t$ the vector
$$
	\xi_t = (x+ty)(1+t^2)^{-1/2}+a
$$
is isotropic. By the condition $\bar\xi_t$ is also isotropic for each sufficiently
small $t$. This implies
$$
	\bar g(\bar x,\bar x)+t^2\bar g(\bar y,\bar y)+2t\bar g(\bar x,\bar y) +
			2\sqrt{1+t^2} \{\bar g(\bar x,\bar a)+t\bar g(\bar y,\bar a) \}
			+(1+t^2)\bar g(\bar a,\bar a)=0 \ .
$$
Hence we derive
$$
	\bar g(\bar x,\bar y)=\bar g(\bar y,\bar a)=\bar g(\bar x,\bar a)=0\,,  \quad 
	\bar g(\bar x,\bar x)=\bar g(\bar y,\bar y)=-\bar g(\bar a,\bar a) \,,
$$
from which the assertion follows easily.

An $n$-dimensional nonflat (pseudo-) Riemannian manifold $M$ is said to be a
$K_n^*$-manifold [12], if it has one of the following properties:

1) $M$ is recurrent, i.e. \,$\nabla R=\alpha\otimes R$, where \,$\alpha \ne 0$;

2) $M$ is symmetric ($\nabla R=0$) and there exists a differential form \,$\alpha\ne 0$, such that
$$
	\sum_{cycl\, X,Y,Z} \alpha(X)R(Y,Z,U,V)=0 \ .
$$
A. W\,a\,l\,k\,e\,r  showed in [12], that $\alpha$ is defined by $\alpha(X)=g(\nabla v,X)$, where 
$v$ is a  function (called recurrence-function) and $\nabla v$ denotes
the gradient of $v$.

An \,$n$-dimensional (pseudo-) Riemannian manifold $M$ is said to be of quasi-constant
curvature [1, 2], if it is conformally flat and there exist functions $H,\,N$ 
and a unit vector $V$, such that the curvature tensor has the form
$$
	R=(N-H)\varphi(B)+H\pi_1 \ ,
$$
where $B(X,Y)=g(X,V)g(Y,V)$. Note that for any point $p$ of $M$ we have $K(\alpha)=H(p)$
for any nondegenerate plane $\alpha$ in $T_pM$, perpendicular to $V_p$ and
$K(\alpha)=N(p)$ for any nondegenerate plane $\alpha$ in $T_pM$, containing $V_p$.
Such a manifold we shall denote by $M(H,N,V)$. Of course, if dim\,$M\ge 4$ and 
the curvature tensor has the above mentioned form, the manifold is necessarily
conformally flat.

Two standart examples of manifolds of quasi-constant curvature are the following
(see also [2]).

E\,x\,a\,m\,p\,l\,e \,1. Let $M_1(c)$ be an $(n-1)$-dimensional (pseudo-) Riemannian
manifold of constant curvature and $M=M_1(c)\times \R$. Then $M$ is of quasi-constant 
curvature.

E\,x\,a\,m\,p\,l\,e \,2. Let $\R_s^{n+1}$ be the pseudo-Euclidean space with an indefinite 
metric of the signature $(s,n+1-s)$. Let $M$ be the $n$-dimensional indefinite
hypersurface given by the equations
$$
	\left\{
		\begin{array}{l}
			x^i=2y^iy^n/\Delta   \qquad i=1,\hdots,n-1  \\
			x^n=y^n(\Delta-2)/\Delta   \\
			x^{n+1}=f(y^n)
		\end{array}
	\right.
$$
for $\Delta>0,\, y^n>0$, where $\Delta=-\sum_{i=1}^s (y^i)^2 +\sum_{i=s+1}^{n-1} (y^i)^2+1$
and $f$ is a smooth function. Then $M$ is a manifold of quasi-constant curvature and 
$$
	H=\frac{f'^2}{(y^n)^2(1+f'^2)} \,, \qquad N=\frac{4f'f''}{y^n(1+f'^2)} \ .
$$

In the following section we shall use the well known fact [5], that if the metrics
$g$ and $\bar g$ on $M$ are related by $\bar g=e^{2\sigma}g$ (thus $(M,g)$ and
$(M,\bar g)$ are conformal), then
$$
	\overline R=e^{2\sigma}\{ R+\varphi(Q) \} \ ,    \leqno (2.8)
$$
where
$$
	Q(X,Y)=X\sigma Y\sigma -g(\nabla_X\nabla\sigma,Y)-\frac12 g(\nabla\sigma,\nabla\sigma)g(X,Y) \ .  
	                  \leqno (2.9)
$$

\vspace{0.2cm}
{\bf 3. Main results.} We begin this section with a theorem, which follows immediately from
Lemmas 2 and 4.

T\,h\,e\,o\,r\,e\,m\, 1. {\it Let $M$ and $\overline M$ be pseudo-Riemannian manifolds of 
dimension $n\ge 3$ and let $M$ be nowhere of constant sectional curvature. If $f$ is 
a diffeomorphism of $M$ onto $\overline M$ satisfying (2.1) for each weakly degenerate plane 
$\alpha_0$ on $M$, then $f$ is conformal.}

If the diffeomorphism $f$ is conformal, we have $f^*\bar g=e^{2\sigma}g$ or
$f^*\bar g=-e^{2\sigma}g$ for a smooth function $\sigma$. Then without loss of generality we
may identify $M$ with $\overline M$ via $f$ and assume $\bar g = e^{2\sigma}g$. We state: 

T\,h\,e\,o\,r\,e\,m\, 2. {\it Under the conditions of Theorem 1 the following propositions hold:

a) if \,$\nabla\sigma$ vanishes identically or $n\ge4$ and $M$ is nowhere conformally flat,
then $f$ is an isometry;

b) if \,$\nabla\sigma$ is isotropic and either $n\ge 4$ or $M$ is conformally flat of dimension 
$n=3$, then $M$ is a conformal flat $K_n^*$-space and $\sigma$ is a function of the 
recurrence-function;

c) if $\left\Arrowvert\nabla\sigma\right\Arrowvert^2=g(\nabla\sigma,\nabla\sigma)$ doesn't vanish
and either $n\ge 4$ or $M$ is conformally flat of dimension $n=3$, then $M$ is a manifold
$M(H,N,V)$ of quasi-constant curvature, where $\nabla H$, $\nabla N$ and $V$ are proportional 
to $\nabla\sigma$.}

P\,r\,o\,o\,f. \,Since $\bar g=e^{2\sigma}g$, (2.1) reduces to
$$
	\overline R(x,\xi,\xi,x)=e^{4\sigma} R(x,\xi,\xi,x)
$$
for arbitrary $p$ in $M$ and $x,\,\xi$ in $T_pM$ with $g(x,x)\ne 0$, $g(\xi,\xi)=g(x,\xi)=0$.
Hence, using Lemma A we obtain
$$
	\overline R=e^{4\sigma} \{ R+c\pi_1 \}
$$
for a function $c$, from which it follows
$$
	\overline R=e^{4\sigma} \{R+(\bar\tau -\tau)\pi_1/(n(n-1)) \}.   \leqno (3.1)
$$

Let $\nabla\sigma=0$, i.e. $\sigma$ is a constant. Then (2.8), (2.9) and (3.1) imply
$(e^{-2\sigma}-1)R=(\bar\tau-\tau)\pi_1/(n(n-1))$. Since $M$ cannot be of constant 
sectional curvature in an open set, this yields $\sigma =0$, i.e. $f$ is an isometry.

On the other hand, if $x,\,y,\,z,\,u$ are arbitrary orthogonal vectors at a point $p$
of $M$, from (2.8) and (3.1) we find
$$
	(e^{2\sigma(p)}-1)R(x,y,z,u)=0 \ .    \leqno (3.2)
$$
If $n\ge 4$ and $M$ is nowhere conformally flat, Lemma B and (3.2) imply $\sigma=0$,
i.e. $f$ is an isometry, proving a).

Let us assume that $\nabla\sigma$ is isotropic or $\left\Arrowvert\nabla\sigma\right\Arrowvert^2$ doesn't
vanish. In both cases if $p$ is an arbitrary point, $\sigma$ cannot vanish in a neighbourhood of
$p$, i.e. there exists a sequence $p_i$ converging to $p$, such that $\sigma(p_i) \ne 0$ for
each $i$. Then (3.2) and Lemma B imply that the Weil conformal curvature tensor of $M$ 
vanishes in $p_i$. By continuity it vanishes in $p$ and hence $M$ is conformally flat, 
if $n\ge 4$. Consequently in cases b) and c) $M$ is conformally flat. Using (2.8) and 
(3.1) it is not difficult to get
$$
	\overline P = e^{2\sigma}P \ ,       \leqno (3.3)
$$
$$
	Q=(e^{2\sigma}-1)P/(n-2)+(\bar\tau e^{2\sigma}-\tau)g/(2n(n-1)) \ ,  \leqno (3.4)
$$
where $P=S-\tau/ng$. Since $(M,g)$ and $(M,\bar g)$ are conformally flat, we have [13]
$$
	\begin{array}{c}
		\ds(\nabla_X(S-\frac{\tau}{2(n-1)}g))(Y,Z)-(\nabla_Y(S-\frac{\tau}{2(n-1)}g))(X,Z)=0 \ ,  \\
		\ds(\overline\nabla_X(\overline S-\frac{\bar\tau}{2(n-1)}\bar g))(Y,Z)-
							(\overline\nabla_Y(\overline S-\frac{\bar\tau}{2(n-1)}\bar g))(X,Z)=0 \ .
	\end{array}          \leqno (3.5)  
$$
Hence, using (3.4) and $\overline\nabla_XY=\nabla_XY+X\sigma Y+Y\sigma X-g(X,Y)\nabla\sigma$
we obtain
$$
	\begin{array}{c}
		\ds  X\sigma P(Y,Z)-Y\sigma P(X,Z)+\frac{n-2}{2n(n-1)} \{ X(\bar\tau-\tau)g(Y,Z)-
					Y(\bar\tau-\tau)g(X,Z) \}  \\
		+ g(X,Z)P(Y,\nabla\sigma)-g(Y,Z)P(X,\nabla\sigma)=0 \ ,
	\end{array}          
$$
which implies immediately
$$
	P(Y,\nabla\sigma)=\frac{n-2}{2n^2} Y(\bar\tau - \tau) \ ,  \leqno (3.6)
$$
$$
	 X\sigma P(Y,Z)-Y\sigma P(X,Z)+\frac{n-2}{2n^2(n-1)} \{ X(\bar\tau-\tau)g(Y,Z)-
					Y(\bar\tau-\tau)g(X,Z) \} =0 \ .      \leqno (3.7)
$$
To prove b) we assume that $\nabla\sigma$ is isotropic and we put in (3.7)
$X=Z=\nabla\sigma$, $Y=\nabla(\bar\tau-\tau)$. Using (3.6) we find
$(\nabla\sigma)(\bar\tau-\tau)=0$. Hence, substituting in (3.7) $X$ by $\nabla\sigma$
we obtain easily $Y(\bar\tau-\tau)Z\sigma=0$ for arbitrary vector fields $Y,\,Z$ on $M$. 
Since $\nabla\sigma$ can not vanish this shows that $\bar\tau-\tau$ is a constant.
Since $\nabla\sigma$ is isotropic, (2.9) implies $Q(X,\nabla\sigma)=0$. Thus,
using (3.4) and (3.6) we derive
$$
	\bar\tau e^{2\sigma}-\tau=0 \ .   \leqno (3.8)
$$
On the other hand, since $\bar\tau-\tau$ is a constant, (3.7) implies
$ X\sigma P(Y,Z)-Y\sigma P(X,Z)=0$ and hence $P(Y,Z)=hY\sigma Z\sigma$ for a smooth
function $h$. Then (2.9), (3.4) and (3.8) yield 
$g(\nabla_X\nabla\sigma,Y)=(1+(1-e^{2\sigma})h/(n-2))X\sigma Y\sigma$. Hence we obtain easily
$$
	(\nabla_XP)(Y,Z)=\{Xh+2h(1+(1-e^{2\sigma})h/(n-2))X\sigma \}Y\sigma Z\sigma \ .
$$
We compare this with (3.5) to find
$$
	\{XhY\sigma-YhX\sigma\}Z\sigma +\frac{n-2}{2n(n-1)} \{X\tau g(Y,Z)-Y\tau g(X,Z) \}=0 \ .
$$
Here we assume, that $X$ is any vector field on $M$, that $Z$  is orthogonal to 
$\nabla\sigma$, $X$ and that $Y$ is not orthogonal to $Z$. It follows that $X\tau=0$, 
i.e. $\tau$ is a constant. Thus, differentiating (3.8) and using $\bar\tau-\tau=\rm const.$ 
we obtain $\bar\tau X\sigma=0$ for any vector field $X$ on $M$. Since $\nabla\sigma$ 
can not vanish this implies $\bar\tau=0$. According to (3.8) $\tau$ vanishes too. 
Then (3.1) shows that $f$ is sectional-curvature-preserving and b) follows from [10].   

Finally, we assume that $\left\Arrowvert\nabla\sigma\right\Arrowvert^2$ doesn't vanish. Let in
(3.7) $X=Z=\nabla\sigma$ and $Y$ be orthogonal to $\nabla\sigma$. Using (3.6) we
obtain $Y(\bar\tau-\tau)=0$. Hence $\mu\nabla\sigma = \nabla(\bar\tau-\tau)$, where
$\mu$ is a smooth function. Substituting $X$ by $\nabla\sigma$ in (3.7) and making
use of (3.6) we find
$$
	P(Y,Z)= \mu \{ \lambda Y\sigma Z\sigma - \frac{n-2}{2n^2(n-1)}g(Y,Z) \} \ ,  \leqno (3.9)
$$
where $\lambda=(n-2)/(2n(n-1)) \left\Arrowvert\nabla\sigma\right\Arrowvert^2$. From (2.9), (3.4)
and (3.9) we derive
$$
	g(\nabla_X\nabla\sigma,Y)=(1+\lambda\mu\frac{1-e^{2\sigma}}{n-2})X\sigma Y\sigma
	                  +\nu g(X,Y) \ ,  \leqno (3.10)
$$
where
$$
	\nu = \mu\frac{e^{2\sigma}-1}{2n^2(n-1)}-\frac12 \left\Arrowvert\nabla\sigma\right\Arrowvert^2
	                       -\frac{\bar\tau e^{2\sigma}-\tau}{2n(n-1)} \ .
$$
In particular, this implies $X\left\Arrowvert\nabla\sigma\right\Arrowvert^2=0$ for any $X$
orthogonal to $\nabla\sigma$. From (3.9) and (3.10) one gets
$$
	(\nabla_XP)(Y,Z)=\frac{X\mu}{\mu}P(Y,Z)+\mu X\lambda Y\sigma Z\sigma+
	       2\lambda\mu(1+\lambda\mu\frac{1-e^{2\sigma}}{n-2})X\sigma Y\sigma Z\sigma
$$
$$
	+\lambda\mu\nu \{ g(X,Y)Z\sigma + g(X,Z)Y\sigma \} \ .
$$
Together with (3.5) this yields
$$
	\begin{array}{c}
		\ds \lambda X\mu Y\sigma Z\sigma +\frac{n-2}{2n(n-1)}(X\tau -\frac1n X\mu)g(Y,Z) \\
		-\ds \{ \frac{n-2}{2n(n-1)}Y\tau - \frac{n-2}{2n^2(n-1)}Y\mu -\lambda\mu\nu Y\sigma\}g(X,Z)=0 \ ,
	\end{array}          \leqno (3.11)
$$
when $\lambda$ is orthogonal to $\nabla\sigma$. Let in (3.11) $g(Y,Y)\ne 0$; $Y$ be orthogonal 
to $X$, $\nabla\sigma$; $Z=\nabla\sigma+Y$. Then we obtain $X\tau=1/nX\mu$. Now we assume
in (3.11) $Y=Z=\nabla\sigma$. The result is $X\mu=0$. Consequently $\nabla\sigma$, $\nabla\tau$
and $\nabla\mu$ are proportional. Since $M$ is conformally flat (3.9) shows that it is a manifold
of quasi-constant curvature $M(H,N,V)$ with $V$ proportional to $\nabla\sigma$. Moreover,
we compute
$$
	H=\frac1{n(n-1)}\{ \tau -\frac1n\mu \} \qquad N=\frac1{n(n-1)}\{\tau+\frac{n-2}{2n}\mu \}.
$$
Consequently $\nabla H$, $\nabla N$ and $\nabla\sigma$ are proportional. This proves the theorem.

R\,e\,m\,a\,r\,k. \,Putting in (3.11) $Z=X$ we obtain
$$
	\mu\nu\nabla\sigma= \left\Arrowvert\nabla\sigma\right\Arrowvert^2(\nabla\tau - \frac1n\nabla\mu) \ .  \leqno(3.12)
$$
If $H-N$ doesn't vanish (3.10) and (3.12) imp[ly
$$
	g(\nabla_XY,\nabla H)=\frac12 \frac1{H-N}\left\Arrowvert\nabla H\right\Arrowvert^2g(X,Y)
$$ 
for all vector fields $X,\,Y$, orthogonal to $\nabla\sigma$. This holds for each manifold of
quasi-constant curvature, see [1] (indeed in [1] this is proved in the definite case, but
there is no principle difference). Note also that in the case $H=$\,const it follows from (3.10) 
and (3.12) that $M$ is locally a product $M_1\times M_2$ in a neighbourhood of any point in
which $M$ is not of constant sectional curvature, where $M_1$ is an $(n-1)$-dimensional
manifold of constant sectional curvature $H$.
 
The following result shows that there exists no nontrivial diffeomorphisms, satisfying (2.1)
for strongly degenerate planes.

T\,h\,e\,o\,r\,e\,m\, 3. {\it Let $M$ and $\overline M$ be pseudo-Riemannian manifolds, such that
$M$ is nowhere conformally flat and of signature $(s,n-s)$ where $s\ge 2,\, n-s\ge 2$. Let $f$ be
a diffeomorphism of $M$ onto $\overline M$, satisfying (2.1) for each strongly degenerate plane
$\alpha_0$ on $M$. Then $f$ is an isometry.}

P\,r\,o\,o\,f. \ According to Lemmas 3 and 4 $f$ is conformal. So without loss of generality we may 
identify $M$ with $\overline M$ via $f$ and assume $\bar g = \varepsilon e^{2\sigma}g$, where
$\varepsilon = \pm1$. Let $\varepsilon =1$ and let $\{\xi,\eta\}$ be an arbitrary orthogonal pair of
isotropic vectors on $M$. Then (2.1) yields
$$
	\overline R(\xi,\eta,\eta,\xi) = e^{4\sigma}R(\xi,\eta,\eta,\xi) \ . \leqno (3.13)
$$ 
On the other hand (2.8) implies
$$
	\overline R(\xi,\eta,\eta,\xi) = e^{2\sigma}R(\xi,\eta,\eta,\xi) \ . \leqno(3.14) 
$$ 
From (3.13) and (3.14) we obtain
$$
 (e^{2\sigma}-1)R(\xi,\eta,\eta,\xi)=0 \ . \leqno (3.15)
$$ 
Since $M$ is nowhere conformally flay, then from (3.15) and by applying Lemma C if follows
that $\sigma=0$. Similar arguments show that the case $\varepsilon=-1$ is not possible. This
completes the proof.

\vspace{0.4in}
\centerline{\large R E F E R E N C E S}

\vspace{0.2in}
\noindent
1. T. A\,d\,a\,t\,i, \,Y. W\,o\,n\,g. Manifolds of quasi-constant curvature I: A manifold
of quasi-

constant curvature and an $s$-manifold. {\it TRU. Math}., {\bf 21}, 1985, 95-103.

\noindent
2. V. B\,o\,j\,u, \,M. P\,o\,p\,e\,s\,c\,u: Espaces \`a courbure quasi-constante. {\it J. Differ. Geom., } {\bf 13},

 1978, 375-383.

\noindent
3. A. B\,o\,r\,i\,s\,o\,v, \,G. G\,a\,n\,c\,h\,e\,v, \,O. Ka\,s\,s\,a\,b\,o\,v. Curvature properties and isotropic 

planes of Riemannian and almost Hermitian manifolds of indefinite metrics.
{\it Ann. Univ. 

Sof., Fac. Math.  M\'ec.,} {\bf 78}, 1984, 121-131.

\noindent
4. M. D\,a\,j\,c\,z\,e\,r, \,K. N\,o\,m\,i\,z\,u. On sectional curvature of indefinite metrics. {\it Math. Ann.,} 

{\bf 247}, 1980, 279-282. 

\noindent
5. L. P. E\,i\,s\,e\,n\,h\,a\,r\,t. Riemannian geometry. Princeton. University Press, 1949.

\noindent
6. S. H\,a\,r\,r\,i\,s. A triangle comparizon theorem for Lorentz manifolds. {\it Indiana Math. J.,}

{\bf 31}, 1982, 289-308.
 
\noindent
7. R. S. K\,u\,l\,k\,a\,r\,n\,i. Curvature and metric. {\it Ann. of Math}., {\bf 91}, 1970, 311-331.

\noindent
8. R. S. K\,u\,l\,k\,a\,r\,n\,i. Curvature structures and conformal transformations. {\it J. Differ. 

Geom}., {\bf 4}, 1970, 425-451.

\noindent
9. R. S. K\,u\,l\,k\,a\,r\,n\,i. Equivalence of Kaehler manifolds and other equivalence problems. 

{\it J. Differ. Geom.,} {\bf 9}, 1974, 401-408. 

\noindent
10. B. R\,u\,h. Krummungstreue Diffeomorphismen Riemannscher und pseudo-Riemannscher

\ Mannigfaltigkieten. {\it Math. Z.,} {\bf 189}, 1985, 371-391.

\noindent
11. J. A. S\,c\,h\,o\,u\,t\,e\,n. Ricci calculus. Berlin, 1954.

\noindent
12. A. G. W\,a\,l\,k\,e\,r. On Ruse's spaces of recurrent curvature. {\it Proc. Lond. Math. Soc., 

\ II S\'er.,} {\bf 52}, 1951, 36-64.

\noindent
13. H. W\,e\,i\,l. Zur Infinitesimalgeometrie. {\it G\"ottingen Nachr.,} {\bf 1921}, 99-112.

\noindent
14. S. T. Y\,a\,u. Curvature preserving diffeomorphisms. {\it Ann Math}., {\bf 100}, 1974, 121-130.

\vspace {0.5cm}
\noindent
{\it Center for mathematics and mechanics \ \ \ \ \ \ \ \ \ \ \ \ \ \ \ \ \ \ \ \ \ \ \ \ \ \ \ \ \ \ \ \ \ \
Received 07.06.1988

\noindent
1090 Sofia   \ \ \ \ \ \ \ \ \ \ \ \ \ \ \ \ \  P. O. Box 373}

\end{document}